\documentstyle [11pt]{article}
\input{mssymb}
{\rm \trivlist \item[\hskip\labelsep{\sc Proof.}]}%
{\hspace*{\fill}$\Box$ \endtrivlist}

\newcommand\bbR{{\Bbb R}}
\newcommand\bbC{{\Bbb C}}

\begin{document} 
    
\title{A compact symmetric symplectic non-Kaehler manifold \\ {\tt
dg-ga/9601012}} 

\author{Eugene Lerman \thanks {Department of Mathematics, Univ.\ of
Illinois, Champaign, IL 61801\ lerman@math.uiuc.edu}}

\maketitle

\begin{abstract}
In this paper I construct, using off the shelf components, a compact
symplectic manifold with a non-trivial Hamiltonian circle action that
admits no Kaehler structure. The non-triviality of the action is
guaranteed by the existence of an isolated fixed point.

The motivation for this work comes from the program of classification
of Hamiltonian group actions.  The Audin-Ahara-Hattori-Karshon
classification of Hamiltonian circle actions on compact symplectic
4-manifolds showed that all of such manifolds are Kaehler. Delzant's
classification of $2n$-dimensional symplectic manifolds with
Hamiltonian action of $n$-dimensional tori showed that all such
manifolds are projective toric varieties, hence Kaehler.  An example
in this paper show that not all compact symplectic manifolds that
admit Hamiltonian torus actions are Kaehler.  Similar technique allows
us to construct a compact symplectic manifold with a Hamiltonian
circle action that admits no invariant complex structures, no
invariant polarizations, etc.
\end{abstract}

\section{Introduction}

The main goals of this paper is to construct a nontrivial example of
a Hamiltonian group action on a compact symplectic manifold which is
not Kaehler.  

The motivation for this work comes from the program of classification
of Hamiltonian group actions.  Delzant proved that if a compact
symplectic manifold admits a completely integrable torus action then
the manifold is a projective toric variety and so, in particular, it
is Kaehler.  Karshon has recently finished a classification of compact
symplectic 4-manifolds admitting a Hamiltonian circle action, the
classification that was initiated by Audin and by Ahara and Hattori.
All such 4-manifolds turn out to be Kaehler.  This should be
contrasted with Gompf's results \cite{G}.  One naturally wonders to
what extent symmetries force symplectic manifolds to be Kaehler.

Ideally a nontrivial example of say a non-Kaehler circle symmetry
would be an example of a Hamiltonian action all of whose fixed points
are isolated.  This would rule out such examples as products of
compact non-Kaehler symplectic manifolds with a two-sphere rotated
about a fixed axis.  Since I started the work on this project, there
has been two interesting developments.  S. Tolman constructed a
compact symplectic 6 manifold which admits a Hamiltonian action of a
2-torus with isolated fixed points and which doesn't admit any
invariant Kaehler structure \cite{T}.  C. Woodward showed that
Tolman's manifold admits a Hamiltonian action of $U(2)$ \cite{W},
thereby showing that symplectic analogues of spherical varieties need
not be Kaehler.  Similar results have been independently obtained by
F. Knop. However it appears that Tolman's manifold does admit a
non-invariant Kaehler structure \cite{T2}.

The example below is a twelve dimensional compact symplectic manifold
with a Hamiltonian circle action which has one isolated fixed point
and two connected fixed submanifolds of dimensions four and ten
respectively. For this manifold $M$, the third Betti number is odd.
Thus $M$ admits no Kaehler structure, invariant or otherwise.

\section{The construction}

McDuff showed that if one embeds symplectically a certain 4 manifold
$N$ into the standard projective 5-space ${\bbC P}^5$ and then blows
up ${\bbC P}^5$ symplectically along $N$, the resulting manifold
$\widetilde{\bbC P}^5_N$ has an odd third Betti number and hence
admits no Kaehler structure.  Consider the action of the circle $S^1$
on ${\bbC P}^6$ given by
$$
\lambda \cdot [z_0, z_1, \ldots, z_6] = [\lambda z_0, z_1, \ldots, z_6].
$$
The action is Hamiltonian, and is free away from the point $p=
[1,0,\ldots 0]$ and the hyperplane $H = \{[0, z_1, \ldots,
z_6]\}\simeq {\bbC P}^5$. Embed the manifold $N$ (symplectically and
equivariantly) into ${\bbC P}^6$ by embedding it into the hyperplane
$H$.  

Guillemin and Sternberg pointed out that symplectic blowing up of
submanifolds can be done equivariantly \cite{GS2}. Blow up ${\bbC
P}^6$ along $N$ equivariantly and symplectically.  The result is a
compact symplectic manifold $M$ with a Hamiltonian circle action.  It
is not hard to see that the fixed point set $M^{S^1}$ has three
components: the point $p$, the manifold $\widetilde{\bbC P}^5_N$ and a
copy of $N$. Since the manifold $N$ admits no Kaehler structure, it
follows immediately that $M$ admits no {\sl invariant} Kaehler
structure.

To show that $M$ admits no Kaehler structure whatsoever, we use Morse
theory.  The moment map of the action of $S^1$ on $M$ is Bott-Morse
and its critical set is precisely the set of fixed points $M^{S^1}$.
It follows from an elementary argument (relative Morse lemma,
Mayer-Vietoris and excision) that $H^3(M) \simeq H^3 (\widetilde{\bbC
P}^5_N) = {\bbR }^3$. \\ 

\subsection*{Conclusions } By a theorem of Gromov and Tischler \cite{Gr},
\cite{Ti}, any compact integral symplectic manifold $X$ of dimension
$n$ can be symplectically embedded in ${\bbC P}^{n+1}$. The above
construction produces a compact symplectic manifold $M_X$ with a
Hamiltonian $S^1$ action that has $X$ as one of the components of the
fixed point set. Thus if the manifold $X$ admits no complex structure
\cite{FGG}, \cite{Geiges} then the manifold $M_X$ admits no {\sl
invariant} complex structure; if the manifold $X$ admits no
polarization \cite{Gotay} then the manifold $M_X$ admits no invariant
polarization; and so on.

\subsection*{Acknowledgments}
The idea for the construction of a symplectic manifold with a
Hamiltonian circle action and no Kaehler structure came to me in a
conversation with Yael Karshon.  Sue Tolman helped me with the
cohomology computations.  

A large part of the work on this project was done at MIT partially
supported by an NSF postdoctoral fellowship.  The author is grateful to
MIT and the NSF for the support.


\begin{thebibliography}{Gompf}

%\bibitem[BL]{BL} L. Bates, E. Lerman, Proper group actions and
%symplectic stratified spaces, {\tt dg-ga/9407003}.

\bibitem[D]{D} T.\ Delzant,  Hamiltoniens p\'{e}riodiques et images convexes de
	l'application moment, {\sl Bull.\ math.\ France} {\bf 116}
	 (1988), 315--339.

\bibitem[FGG]{FGG} M. Fern\'andes, M. J. Gotay, and A. Gray, Compact
parallelizable four-dimensional symplectic and complex manifolds, {\sl
Proc.\ A.M.S.} {\bf 103} (1988), 1209--1212.

\bibitem[Ge]{Geiges} H. Geiges, Symplectic structures on $T^2$ bundles
over $T^2$, {\sl Duke Math.\ J.} {\bf 67} (1992), 539--555.

\bibitem[Gom]{G} R.E. Gompf, A new construction of symplectic manifolds,
              preprint, {\sl Ann.\ Math.} {\bf 142} (1995), 527--598.

\bibitem[Got]{Gotay} M. J. Gotay, A class of non-polarizable
symplectic manifolds, {\sl Monatsch.\ Math.} {\bf 103} (1987), 27--30.

\bibitem[Gr]{Gr} M. Gromov, {\sl Partial Differential Relations},
Springer Verlag; New York, Berlin, 1986.
            
%\bibitem[GS1]{GS1} V.\ Guillemin and S.\ Sternberg, Geometric quantization 
%	and multiplicities of group representations, {\sl Invent.\ Math.}
%        {\bf 67} (1982), 515--538.

\bibitem[GS]{GS2} V. Guillemin and S.\ Sternberg,  Birational equivalence 
	in the symplectic category, {\sl Invent.\ Math.} {\bf 97} (1989), 
	485--522.

%\bibitem[K]{K} K. Kodaira, On the structure of compact complex
%analytic surfaces I, {\sl Amer.\ J. Math.} {\bf 86} (1964), 751--798.

%\bibitem[L]{L} E. Lerman,   Symplectic cuts, {\sl Math.\ Research Lett.}
%(1995), {\bf 2}, 247--258.

%\bibitem[Ma2]{Marle} C.-M. Marle, Sous-veri\'et\'e de rang constant
%d'une veri\'et\'e symplectiques, {\sl Ast\'erisque } {\bf 107--108}
%(1983), 69--86.


\bibitem[McD]{McD} D. McDuff, Examples of simply connected non Kaehlerian 
manifolds, {\sl J. Diff.\ Geom.} {\bf 20} (1984), 267--277. 

\bibitem[Ti]{Ti} D. Tischler, Closed 2-forms and an embedding theorem
for symplectic manifolds, {\sl J. Diff.\ Geom.} {\bf 12} (1977),
229--235.

\bibitem[T]{T} S. Tolman, Examples of non-Kaehler Hamiltonian torus
actions,  {\tt dg-ga/9511007, URL http://xxx.lanl.gov/abs/dg-ga/9511007}.

\bibitem[T2]{T2} S. Tolman, personal communication, November 1995.
 
\bibitem[W]{W} C. Woodward, Multiplicity-free Hamiltonian actions need
not be K\"ahler, {\tt dg-ga/9506009, URL
http://xxx.lanl.gov/abs/dg-ga/9506009}.

 

\end{thebibliography}
\end{document}